\documentclass[leqno,12pt]{article}
\newcommand{\fin}{\hfill \framebox[2.4mm]{ }\rule[-3mm]{0mm}{1cm}}
\newcommand{\N}{{\rm I\!N}}
\newcommand{\R}{{\rm I\!R}}
\newcommand{\E}{{\rm I\!E}}
\newcommand{\PP}{{\rm I\!P}}
\newcommand{\D}{{\rm I\!D}}

\newcommand{\dd}{{\rm d\!I}}
\newtheorem{proposition}{Proposition}
\newtheorem{lemme}[proposition]{Lemma}
\newtheorem{theoreme}[proposition]{Theorem}

\newtheorem{corollaire}[proposition]{Corollary}
\title{On numerical integration by the shift and application to  Wiener
space} \author{Nicolas BOULEAU \thanks{CERMA, Ecole Nationale des Ponts et
Chauss\'{e}es 93167 Noisy-le-Grand cedex France}} 
\date{-=-}

\begin{document}
\maketitle

\vspace{2cm}

Since the advantages of quasi-Monte Carlo methods vanish when the dimension of
the basic space increases, the question arises whether there are better methods
than classical Monte Carlo in large or infinite dimensional basic
spaces. We study here the use of the shift operator with the pointwise ergodic
theorem whose implementation is particularly interesting. After recalling the
theoretical results on the speed of convergence in a form useful for
applications, we give sufficient criteria for the law of iterated logarithm in
several cases and in particular in situations involving the  Wiener space.

If a family of real random variables is naturally defined on a probability
space which can be smoothly changed to be $([0,1]^s,{\cal B}([0,1]^s),dx)$
with small $s$, quasi-Monte Carlo methods are among the fastest ones for computing expectations, at
least when the family is wide enough to exclude other specific methods.
See for example \cite{niederreiter1},
\cite{niederreiter2}. But the advantage of these methods vanishes
when $s$ increases (cf \cite{sarkar}). Practically, for the best
low-discrepancy sequences available at present (cf
\cite{niederreiter2}, \cite{h.faure}), to compute expectations with an
accuracy of $10^{-4}$ by unit of standard deviation, it is faster to
come back to the classical Monte Carlo method as soon as the dimension
$s$ exceeds $20$ (cf \cite{benalaya} \cite{bouleau-pages}).

In large or infinite dimension (computation of expectations of stopping
times for Markov chains, or of functionals of solutions of stochastic
differential equations, etc.) the classical Monte Carlo method which is
based on the law of large numbers, can nevertheless be supplanted by another method based on the
pointwise ergodic theorem of Birkhoff and the shift operator. Particular features of the implementation of
this method make it at present the most interesting way of integration in large or
infinite dimension (cf \cite{bouleau}).

The aim of this study is to clarify the consequences of recent theoretical
results for the numerical computation of expectation by the shift
method, and in particular  to yield sufficient criteria for the existence of
speed of convergence of the type `iterated logarithm' in several 
situations. We put particular accent to the case of Wiener space because
it is the basic space of many situations useful in  applications. 

The content of the study is the following:

\vspace{0.5 ex}I. Law of iterated logarithm for the shift

II. Criteria of membership for the Gordin class

\hspace{0.5cm} 1. Case of the torus $T^s$

\hspace{0.5cm} 2. Case of the torus $T^{\N}$

\hspace{0.5cm} 3. Case of Wiener space

\hspace{1 cm} a) The Wiener space as a product space

\hspace{1 cm} b) Functionals of lipschitzian SDE's

\hspace{1 cm} c) Multiple Wiener integrals

\hspace{0.5 cm} 4. Other factorisations of the Wiener space

\vspace{0.3 ex}We give now some details on each of these  parts.

The  first part is concerned with the speed of convergence in the
pointwise ergodic theorem for the shift on $T^{\N}$. In contrast to the
case of the law of large numbers, there is no standard speed of
convergence valid for every function in $L^2$. Nevertheless the
successive improvements of the LIL (cf \cite{heyde-scott},
\cite{scott}, \cite{hall-heyde}, and more recently \cite{berger}) have
shown the importance of a sub-class of $L^2$ for which a form of the LIL
is valid and which contains several useful examples (cf part II). We call
this class the Gordin class by reference to \cite{gordin} one of
the first works where this  decomposition in sum of martingale
increments and a subsidiary harmless term is used. Our purpose is not to
extend the  general results (cf \cite{berger}) but   to
express useful consequences for applications. All the results are
explicitely proved except the theorem of Heyde and Scott itself.

In the  second part, we show first that functions in the Sobolev
spaces $H^\alpha(T^s)$ are in the Gordin class for the shift of binary
digits. Next for the torus $T^{\N}$ with the shift to the right, Dirichlet
forms techniques are used to obtain a simple sufficient criterion for
membership to the Gordin class. For the first factorisation of the Wiener
space under study, the shift becomes the scaling $B_t\circ \tau
=\frac{1}{\sqrt{2}}B_{2t}$. With this transform, H\"{o}lderian functions
of solutions of Lipschitzian SDE's are shown to belong to the Gordin
class. Somes examples are analysed which are related to multiple
Wiener integrals. Other  factorisations are discussed and especially the
representation of  Brownian motion  on the Schauder basis of ${\cal
C}([0,1])$ consisting of primitive functions of the Haar basis. The
criterion obtained on $T^{\N}$ by Dirichlet forms method applies here as
well.

We thank J.P. Conze and E. Lesigne for useful discussions and
suggestions. 
\section{The law of iterated logarithm for the shift }
We are interested by almost sure results. It is well known (see 
\cite{halasz},  \cite{krengel1}) that for every ergodic endomorphism $\tau$
on a Lebesgue space, and for every sequence, $(\alpha_n)$, $\alpha_n>0$,
$\alpha_n\rightarrow 0$, there is an $f$ in $L^2$ such that
$$(\frac{1}{n}\sum_{k=0}^{n-1}f\circ\tau^k-\E f)/\alpha_n \rightarrow +\infty
\hspace{1cm}{\mbox{a.s.}}$$ Such a ``slow" $f$ is constructed by suitable
application of the Rohlin-Halmos lemma which likewise furnishes a ``fast"
non constant $f\in L^2$ for which $(\frac{1}{n}\sum_{k=0}^{n-1}f\circ\tau^k-\E
f)$ is $o(\frac{1}{n^{1-\epsilon}})$ (see also \cite{krengel2} pp.14-15)

Nevertheless such functions, by the nature of the Rohlin-Halmos lemma itself,
are rather abstract examples, and do not prohibit an LIL from holding for a
large class of functions containing the most common ones.

For later convenience, we assume the following framework: 

$(E, {\cal E}, \mu)$ is a probability space and 

$(\Omega, {\cal A}, \PP) = (E^{\bf Z}, {\cal E}^{\otimes{\bf Z}},
\mu^{\otimes{\bf Z}})$.

The coordinates from $\Omega$ into $E$ are denoted by $X_n$. We define the
ergodic automorphism $\tau$ on $(\Omega, {\cal A}, \PP)$ by $$ X_n\circ
\tau=X_{n-1}\hspace{1cm}\forall n\in{\bf Z}.$$

We call $\tau$ the shift to the right. One puts $${\cal
F}^n_m=\sigma(X_m,X_{m+1},\ldots,X_n)\hspace{1cm}{\mbox{for }}m\leq n\in{\bf
Z}$$
$${\cal F}^n_{-\infty}=\sigma(X_k, k\leq n)$$
$${\cal F}^{+\infty}_{m}=\sigma(X_k, k\geq m)$$
$${\cal F}^{+\infty}_{-\infty}=\sigma(X_k, k\in {\bf Z})$$
As stated in the introduction the following results can be proved in a more
general setting, for other endomorphisms (see \cite{heyde-scott},
\cite{scott}, \cite{hall-heyde}) and for Banach-valued random variables
(see \cite{berger}).

Let us consider on $L^1({\cal F}_0^\infty,\PP)$ the Perron-Frobenius operator
$T$ defined by 
$$Tf= \E [f|{\cal F}^\infty_1]\circ\tau\hspace{1cm}f\in L^1({\cal
F}_0^\infty,\PP)$$
we then have:
\begin{lemme}
For $f\in L^2({\cal F}_0^\infty,\PP)$, $\E f=0$, the following assumptions
are equivalent:

a) $\sum_{n=0}^NT^nf$ remains bounded in $L^2$,

b) $\sum_{n=0}^NT^nf$ converges weakly in  $L^2$ when $N\uparrow \infty$,

c) $\sum_{n=0}^NT^nf$ converges  in  $L^2$ when $N\uparrow \infty$,

d) there exists $g\in L^2({\cal F}_0^\infty,\PP)$ such that $f=(I-T)g$.
\end{lemme}

\noindent{\bf Proof}. b)$\Rightarrow$d): if $\sum_{n=0}^NT^nf$ converges
weakly, by the Banach-Steinhaus theorem the limit $g$ is an element of
$L^2({\cal F}_0^\infty,\PP)$. By composition with the bounded operator $T$, we
obtain $g=f+Tg$.

 d)$\Rightarrow$c): if $f=(I-T)g$, $g\in L^2$, it can be supposed $\E g= 0$.
Then $\| T^Ng\|_{L^2}\rightarrow 0$ when $N\uparrow\infty$. Indeed $\|
T^Ng\|_{L^2}^2=\E[\E(g|{\cal F}_N^\infty)^2]$ and $\E(g|{\cal F}_N^\infty)$ is
an inverse martingale which tends to zero in $L^2$.

Finally for a)$\Rightarrow$b), let us consider a subsequence $N_k$ such that 
$\sum_{n=0}^{N_k}T^nf$ converges weakly in $L^2({\cal F}_0^\infty)$ as
$k\uparrow\infty$. Letting $g$ be the limit, by composition with $T$ we get
$$g-f+\lim_{k\uparrow\infty}T^{N_k+1}g=Tg$$ and the same argument as for
c)$\Rightarrow$d) shows that $T^{N_k+1}f\rightarrow 0$ in $L^2$.\fin

\noindent {\bf Remark}. It is easy to see that these conditions are
equivalent to the condition that $\sum_{n=0}^N\tau^nf$ converge for the
topology $\sigma(L^2({\cal F}^{+\infty}_{-\infty}),L^2({\cal
F}_0^{+\infty}))$.

We shall say that a function $f\in L^2({\cal F}_0^{+\infty})$ belongs to
{\bf the Gordin class} (for which we  write $f\in{\cal G}$) if $f-\E f$
satisfies the equivalent conditions of lemma 1.
\begin{lemme}The Gordin class is the class of the functions $f\in L^2({\cal
F}_0^{+\infty})$ admitting a decomposition 
\begin{equation}\label{1}
f-\E f=\tilde{g}+h\circ\tau^{-1}-h
\end{equation}
where $\tilde{g},h\in f\in L^2({\cal F}_0^{+\infty})$ with
$\E(\tilde{g}|{\cal F}_1^\infty)=0$ and $\E h=0$. Such a decomposition, if it
exists, is unique.
\end{lemme}

\noindent {\bf Proof}. By lemma 1, if $f\in{\cal G}$ there is a $g\in 
L^2({\cal F}_0^{+\infty})$ such that $f-\E f=g-Tg$. Putting
$\tilde{g}=g-\E(g|{\cal F}_1^\infty)$ and $h=Tg=\E(g|{\cal
F}_1^\infty)\circ\tau$, we get the decomposition (\ref{1}). 

Conversely, if $f$ can be decomposed as (\ref{1}), we
have $T\tilde{g}=0$ and $T(h\circ\tau^{-1})= h$, hence $f-\E f=(I-T)g$ with
$g=\tilde{g} + h\circ\tau^{-1}$. The uniqueness follows immediately. \fin

The theorem of iterated logarithm is valid for functions in the Gordin class:
\begin{theoreme} Let $f\in L^2({\cal F}_0^{+\infty})$ be in the Gordin
class, and $\tilde{g},\,\,h$ the elements of its decomposition
(\ref{1}). Then, putting $S_N=\sum_{n=0}^N(f-\E
f)\circ\tau^n$, there holds

a)
$$\lim_{N\uparrow\infty}\frac{1}{\sqrt{N}}\|S_N\|_{L^2}=\|\tilde{g}\|_{L^2}$$

b) $$\limsup_{N\uparrow\infty}\frac{|S_N|}{\sqrt{2N\log\log
N}}=\|\tilde{g}\|_{L^2}$$ 
\end{theoreme}

\noindent{\bf Proof}. Noting that 
$$S_N =\sum_{n=0}^N\tilde{g}\circ\tau^n+h\circ\tau^{-1}-h\circ\tau^N$$
part a) comes from the following inequality, where the norms are taken in
$L^2$:
$$\left|\|\frac{1}{\sqrt{N}}S_N\|-\|\frac{1}{\sqrt{N}}\sum_{n=0}^N\tilde{g}\circ\tau^n\|\right|\leq
2\frac{\|h\|}{\sqrt{N}}\rightarrow_{N\uparrow\infty}0$$
and from $\|\sum_{n=0}^N\tilde{g}\circ\tau^n\|^2=(N+1)\|\tilde{g}\|^2$, which
follows by orthogonality.

If $\tilde{g}=0$ part b) is a consequence of the fact that, $h$ being in
$L^2$, $\frac{h\circ\tau^N}{\sqrt{N}}\rightarrow 0$ when $N\uparrow\infty$ by
the pointwise ergodic theorem. Thus, when $\tilde{g}\neq 0$ it suffices to show
that 

$$\limsup_{N\uparrow\infty}\frac{|\sum_{n=0}^N\tilde{g}\circ\tau^n|}{\sqrt{2N\log\log
N}}=\|\tilde{g}\|.$$ But this is given by the theorem of Heyde and Scott
(\cite{heyde-scott} corollary 2). \fin

We shall now state sufficient conditions for membership to the Gordin class
${\cal G}$.

Without subscript, norms are $L^2$-norms.
\begin{proposition} Let $f\in L^2({\cal F}_0^{+\infty})$ be such that 
\begin{equation}\label{2}
\sum_{n=0}^\infty\|\E[f]-\E(f|{\cal F}_n^\infty)\|<+\infty,
\end{equation}
then $f\in{\cal G}$ and the $\tilde{g}$ of its decomposition satisfies 
$$\|\tilde{g}\|\leq \sum_{n=0}^\infty\|\E[f]-\E(f|{\cal F}_n^\infty)\|.$$
\end{proposition}
\noindent{\bf Proof}. By the fact that 
$$\|T^n(f-\E f)\|= \|\E f-\E(f|{\cal F}_n^\infty)\|$$
the convergence of the series (\ref{2}) implies that the series $\sum T^n(f-\E
f)$ converges normally. Letting $g$ be its sum, then $\tilde{g}$ is given by
$g-\E(g|{\cal F}_1^\infty)$ thus $\|\tilde{g}\|\leq\|g\|$.\fin 
\begin{proposition} Let $f\in L^2({\cal F}_0^{+\infty})$, and  let us consider
the decomposition 
\begin{equation}\label{3}
f=\E f +\sum_{n=0}^\infty f_n
\end{equation}
with 
$$f_0=\E (f|{\cal F}_0^0)-\E(f)$$
and 
$$f_n=\E(f|{\cal F}_0^n)-\E(f|{\cal
F}_0^{n-1})\hspace{1cm}{\mbox{for }}n\geq 1$$

a) $f\in{\cal G}$ if and only if 
$$\sup_N\sum_{j=o}^\infty\|\E(\sum_{n=0}^N f_{n+j}\circ\tau^n|{\cal
F}_0^\infty)\|^2<+\infty.$$

b) This is satisfied if 
$$\sum_{m\geq 0}\sqrt{\sum_{k\geq m}\|f_k\|^2}<+\infty$$
and then the $\tilde{g}$ associated with $f$ in (\ref{1}) is such that 
$\|\tilde{g}\|\leq \sum_{m\geq 0}\sqrt{\sum_{k\geq m}\|f_k\|^2}.$

c) This is also satified if 
$$\sum_{m\geq 0}\sqrt{m}\|f_m\|<+\infty$$
and then the $\tilde{g}$ associated with $f$ in (\ref{1}) is such that
$\|\tilde{g}\|\leq \sum_{m\geq 0}\sqrt{m}\|f_m\|$.
\end{proposition}
\noindent{\bf Proof}. The existence of the decomposition (\ref{3}) for any 
$f\in L^2({\cal F}_0^{+\infty})$ comes from the fact that 
$$\sum_{n=1}^Nf_n=\E(f|{\cal F}_0^N)-\E(f|{\cal F}_0^0)$$
is a martingale which converges in $L^2$.

Let $f\in L^2({\cal F}_0^{+\infty})$, and put $\tilde{f}=f-\E f$. Using
the fact that for $n\geq 0$, $T^n \tilde{f}=\E[\tilde{f}\circ\tau^n|{\cal
F}_0^\infty]$
 we get 
$$T^n\tilde{f}=\sum_{k\geq n}\E[f_k\circ \tau^n|{\cal F}_0^\infty]$$
because for $k<n$, $f_k\circ\tau^n$is ${\cal F}_{-\infty}^{-1}$-measurable.
It follows that
 \begin{equation}\label{5}
\sum_{n=0}^N T^n\tilde{f}=\sum_{n=0}^N\sum_{k\geq n}\E[f_k\circ\tau^n|{\cal
F}_0^\infty]
\end{equation}
$$=\sum_{j=0}^\infty\E(\sum_{n=0}^Nf_{n+j}\circ\tau^n|{\cal F}_0^\infty).$$
But the random variables 
$$Z_j^N=\E(\sum_{n=0}^Nf_{n+j}\circ\tau^n|{\cal F}_0^\infty)$$
form an orthogonal sequence and therefore
$$\|\sum_{n=0}^N T^n\tilde{f}\|^2=\sum_{j=0}^\infty\|Z_j^N\|^2$$
and part a) follows from lemma 1.

From the equality (\ref{5}) we have also 
$$\|\sum_{n=0}^NT^n\tilde{f}\|\leq \sum_{n=0}^N\|\sum_{k\geq
n}\E[f_k\circ\tau^n|{\cal F}_0^\infty]\|.$$
For every fixed $n$ the sequence $(\E(f_k\circ\tau^n\|{\cal
F}_0^\infty))_{k\geq n}$ is orthogonal, and so
$$\|\sum_{k\geq
n}\E[f_k\circ\tau^n|{\cal F}_0^\infty]\|=\left(\sum_{k\geq n}\|\E[f_k\circ\tau^n|{\cal
F}_0^\infty]\|^2\right)^{1/2}$$
$$\leq \left(\sum_{k\geq n}\|f_k\|^2\right)^{1/2}$$
which gives part b).

Taking once more the equality (4) rewritten as 
$$\sum_{n=0}^N
T^n\tilde{f}=\sum_{m=0}^\infty\E(f_m+f_m\circ\tau+\cdots+f_m\circ\tau^{m\wedge
N}|{\cal F}_0^N)$$
gives 
$$\|\sum_{n=0}^N
T^n\tilde{f}\|\leq
\sum_{m=0}^\infty\|f_m+f_m\circ\tau+\cdots+f_m\circ\tau^{m\wedge N}\|$$
$$\leq \sum_{m=0}^\infty\sqrt{m}\|f_m\|,$$
by the fact that the sequence $f_m\circ\tau^{m\wedge
N},\;f_m\circ\tau^{m\wedge N-1},\ldots,f_m$ is orthogonal. Part c) follows by
the same arguments.\fin

If $T$ is an almost finite stopping time of the $\sigma$-fields $({\cal
F}_0^n)_{n\geq 0}$ and if $f$ is an ${\cal F}_0^T$-measurable random variable, $f$
can be written as $$f=\sum_{k\geq 0}f\;1_{\{T=k\}}$$
with $f\;1_{\{T=k\}}$ ${\cal F}_0^k$-measurable. This is a particular case of
the following situation:
\begin{proposition} Let $f\in L^2({\cal F}_0^\infty)$ admit the following
representation converging in $L^2$:
$$f=\sum_{k=0}^\infty f_k \hspace{0.5cm}{\mbox{with }}f_k\;\;{\cal
F}_0^k{\mbox{-measurable}}.$$
If the condition 
\begin{equation}
\label{6} \sum_{k=0}^\infty k\|f_k-\E f_k\|<+\infty
\end{equation}
is fulfilled, then $f\in{\cal G}$, and the associated $\tilde{g}$ satisfies
\begin{equation}
\label{7} \|\tilde{g}\|\leq \sum_{k=0}^\infty\sqrt{k+1}\|f_k-\E f_k\|.
\end{equation}
\end{proposition}
\noindent{\bf Proof}. By the fact that $T^n(f_k-\E f_k)=0$ for $n>k$, 
$$T^n(f-\E f)=\sum_{k\geq n}T^n(f_k-\E f_k)$$
and therefore under condition
(5) the series $\sum_n T^n(f-\E f)$ is normally convergent and $f\in
{\cal G}$.

Let us put $g(f_k)=\sum_{n\geq 0}T^n(f_k-\E f_k)$ and $g(f)=\sum_{n\geq
0}T^n(f-\E f)$. Under condition (5) we have thus $$g(f)=\sum_{k\geq 0}
g(f_k),$$ the series converging normally. Therefore
$$g(f)-\E(g(f)|{\cal F}_0^1)=\sum_{k\geq 0}[g(f_k)-\E(g(f_k)|{\cal F}_0^1)],$$
the series again converging normally. But by lemma 7 below and proposition 5
applied to $f_k$ we have 
$$ \|g(f_k)-\E(g(f_k)|{\cal F}_0^1)\|\leq \sqrt{k+1}\|f_k-\E f_k\|$$
so that 
$$\|g(f)-\E(g(f)|{\cal F}_0^1)\|\leq \sum_k\sqrt{k+1}\|f_k-\E f_k\|$$
which proves the proposition.\fin
\begin{lemme} If $f\in L^2$ depends only on $d$ consecutive coordinates then
$$\limsup_N\frac{|f+f\circ\tau+\cdots+f\circ\tau^{N-1}-N\E
f|}{\sqrt{2N\log\log N}}\leq \sqrt{d}\|f-\E f\|$$
\end{lemme}
\noindent{\bf Proof}. This is a simple application of the LIL of
Hartman-Wintner for independent variables. Let us put $N-1=pd+q$ with $0\leq
q<d$, and let us suppose $f$ to be centred. Then 
\begin{equation}
\label{8}
\sum_{i=0}^{N-1}f\circ\tau^i=\sum_{k=0}^{d-1}\sum_{j=0}^{p-1}f\circ\tau^{jd+k}+\sum_{n=pd}^{pd+q}f\circ\tau^n.
\end{equation}
By the fact that for every fixed $k$ 
$$\limsup_p\frac{|\sum_{j=0}^{p-1}f\circ\tau^{jd+k}|}{\sqrt{2p\log\log
p}}=\|f\|$$
we have
 $$\begin{array}{rcl}
\limsup_N\frac{|\sum_{k=0}^{d-1}\sum_{j=0}^{p-1}f\circ\tau^{jd+k}|}{\sqrt{2N\log\log
N}}&\leq & d\|f\|\lim\frac{\sqrt{2p\log\log p}}{\sqrt{2N\log\log N}}\\
&\leq &\sqrt{d}\|f\|.
\end{array}$$
Now the second term of (7) gives 
$$\frac{|\sum_{n=pd}^{pd+q}f\circ\tau^n|}{\sqrt{2N\log\log N}}\leq
\frac{\sum_{j=0}^{d-1}|f|\circ\tau^{pd+j}}{\sqrt{2N\log\log N}}$$
which vanishes almost surely as $N\uparrow\infty$ by the ergodic theorem
because $f\in L^2$. The lemma follows from these estimates.\fin

\noindent{\bf Remark}. For $f\in L^2({\cal F}_0^\infty)$
admitting the representation $f=\sum_k f_k$ converging in $L^2$ with $f_k$
${\cal F}_0^k$-measurable, we dont know whether the sole hypothesis
$\sum_k\sqrt{k+1}\|f_k-\E f_k\|<+\infty$
suffices to imply 
$$\limsup_N\frac{|\sum_{n=0}^N(f-\E f)\circ\tau^n|}{\sqrt{2N\log\log N}}\leq \sum_k\sqrt{k+1}\|f_k-\E
f_k\|.$$

The following result, whose statement is simple, is a rather rough consequence
of the preceding proposition.
\begin{proposition} Let $T$ be an a.s. finite stopping time of $({\cal
F}_0^n)_{n\geq 0}$, and $f\in L^2({\cal F}_0^T).$

If there is an $\alpha>1$ such that 
$$\E[f^2T^3\log^\alpha T]<+\infty$$
then 
$$\limsup_N\frac{|\sum_{n=0}^N(f-\E f)\circ\tau^n|}{\sqrt{2N\log\log N}}\leq
\frac{\sqrt{6}}{\pi}\|f(T+1)^{3/2}\|.$$
 \end{proposition}
\noindent{\bf Proof}. Putting $f_k=f\;1_{\{T=k\}}$ we have 
$$\sum_{k\leq 1} k\|f_k-\E f_k\|\leq \sum_k
k\|f_k\|=\sum_k\frac{1}{\sqrt{k\log^\alpha k}}\|\sqrt{k^3\log^\alpha k}f_k\|$$
which by the  Cauchy-Schwarz inequality, is bounded by 
$$\leq \sqrt{\sum_k\frac{1}{k\log^\alpha k}}\sqrt{\sum_k k^3\log^\alpha k
f_k}\leq c\sqrt{\E[T^3\log^\alpha Tf^2]}.$$
Therefore the inequality (5) is satisfied and similarly 
$$\sum_{k=0}^\infty\sqrt{k+1}\|f_k-\E f_k\|\leq
\frac{\sqrt{6}}{\pi}\|f(T+1)^{3/2}\|.$$\fin

\noindent{\bf Remark}. It is worth noting that if $f\in{\cal G}$ it can of
course happen that 
$$\|\tilde{g}\|<\|f-\E f\|.$$
This is the case if $f-\E f=h\circ\tau^{-1} -h$ with $h\in L^2$ and other
examples are easily constructed by the Gordin decomposition. This can occur
even when $f$ depends only on a finite number of coordinates.
In this case  integration by the shift method runs (asymptotically)
faster than by classical Monte Carlo.

Nevertheless, the  principal interest of the shift method does not come
from this phenomenon but from certain facilities afforded by its implementation
(see \cite{bouleau}). 

\section{Criteria of membership to the Gordin class}
\subsection{The case of the torus $T^s$}
Let us consider the following transform of $T^s\times T^s$
$$((x_1,\ldots,x_s),(y_1,\ldots,y_s))\stackrel{\tau}{\longrightarrow}
(([2x_1],\ldots,[2x_s])),(\frac{2x_1-[2x_1]+y_1}{2},\ldots,(\frac{2x_s-[2x_s]+y_s}{2}))$$
where $[x]$ is the fractional part of $x\in\R$, which is easily seen to
correspond to the bilateral Bernoulli shift by binary expansion of real numbers.

We have for this transformation and for $f\in L^2(T^s, dy_1\cdots dy_s)$: 
$$Tf(y)=\sum_{n\in\{0,1\}^s}\frac{1}{2^s}f(\frac{n}{2}+\frac{y}{2})$$
$$T^nf(y)=\sum_{k\in\{0,\ldots,2^n-1\}^s}\frac{1}{2^{ns}}f(\frac{k}{2^n}+\frac{y}{2^n}).$$
Using then the Fourier representation of $f$ 
$$f(y)= \sum_{m\in{\bf Z}^s}a_m\;e^{2i\pi<m,y>}$$
one easily obtains that if $\E f=0$
$$\sum_{n=0}^N T^nf(y)=\sum_{q\in {\bf Z}^s,\;q\neq 0}\;\;\sum_{n=0}^Na_{2^n
q}\;e^{2i\pi<q,y>}$$
and we get:
\begin{proposition} The function $f\in L^2(T^s)$ belongs to ${\cal G}$ if and
only if 
$$\sup_N \sum_{q\in {\bf Z}^s,\;q\neq 0}\;\;|\sum_{n=0}^Na_{2^n q}|^2<+\infty$$
in which case $\forall q\neq 0$, $\sum_{n=0}^Na_{2^n
q}\rightarrow_{N\uparrow\infty}b_q$ with $\sum_q|b_q|^2<+\infty$
and the $\tilde{g}$ of $f$ in (1) satisfies $\|\tilde{g}\|^2\leq \sum_q |b_q|^2$. 
\end{proposition}
\begin{corollaire} Let $f\in L^2(T^s)$ be such that there are $c_n\geq 0$ with
$\sum_{n=0}^\infty c_n<+\infty$ and $|a_{2^nm}|\leq c_n|a_m|$ $\forall m\in{\bf
Z}^s\backslash\{0\}$. Then $f\in{\cal G}$ and $$\|\tilde{g}\|\leq \|f-\E
f\|\sum_{n=0}^\infty c_n.$$
\end{corollaire}

\noindent{\bf Example}. Letting $f$ belong to the Sobolev space $H^\alpha(T^s)$
defined by $\sum_{p\in{\bf Z}^s}|a_p|^2(\sum_{i=0}^s p_i^2)^\alpha <+\infty$
for some $\alpha>0$. Then $f\in {\cal G}$ and $\|\tilde{g}\|\leq
c(\sum_{p}|a_p|^2(\sum_{i=0}^s p_i^2)^\alpha)^{1/2}.$

\subsection{Case of the infinite dimensional torus}
We consider here the Bernoulli shift ( to the right) on $([0,1]^s,{\cal
B}([0,1]^s),dx)^{\bf Z}.$

The property of membership to the Gordin class is strongly related to the
dependence of $f$ on the size of the derivatives of $f$ (when they exist) with
respect to the faraway coordinates. This is  particularly  simple to  espress
by means of Dirichlet forms:

Let us consider a Dirichlet form $(\dd,\varepsilon)$ on
$L^2([0,1]^s,dx)$ possessing a carr\'{e} du champ operator $\gamma$ (cf
\cite{bouleau-hirsch1}) and let us consider the product Dirichlet structure
(cf \cite{bouleau-hirsch2}):
$$(\Omega,{\cal A},\PP,\dd,{\cal E})=([0,1]^s,{\cal
B}([0,1]^s),dx,\dd,\varepsilon)^{\otimes{\bf Z}}$$
This structure has a carr\'{e} du champ $\Gamma$ given by
$$\Gamma(f,f)=\sum_{i\in{\bf Z}}\gamma_i(f,f) \hspace{0.5cm}\forall f\in\D$$
where $\gamma_i$ operates on the $i$-th coordinate. We consider the shift $\tau$
given by 
$$X_n\circ\tau=X_{n-1}$$
where $(X_n)_{n\in{\bf Z}}$ are the coordinates.
For $F\in L^2(\Omega,{\cal F}_0^\infty)$ there holds 
$$Tf(x_0,x_1,\ldots,x_n,\ldots)=\int_{x\in[0,1]^s}f(x,x_0,x_1,\ldots)dx.$$
We make the following assumption (\ref{**}): Let $L_0^2=\{f\in L^2\;\E f=0\}$
\begin{equation}
\label{**}
\left\{
\begin{array}{l}
\mbox{There exists }K>0\mbox{ such that}\\
\forall f\in \D\cap L_0^2\hspace{0.5cm}\|f\|_{L^2}^2\leq K{\cal E}(f,f)\\
\end{array}
\right. 
\end{equation}
 
Then the space $\D\cap L_0^2$ is  a Hilbert space for the norm $\sqrt{{\cal
E}(f,f)}$ which is invariant by $\tau$.
Let $\D_0=\{f\in\D,\;\E f=0,\; f \mbox{ is }{\cal
F}_0^\infty{\mbox{-measurable}}\}$
which is closed in $\D$
\begin{proposition} Under hypothesis (\ref{**}), let $f\in\D_0$ be such that 
$$\sum_{k=0}^\infty(\sum_{i=k}^\infty\E [\gamma_i(f,f)])^{1/2}<+\infty.$$
Then $f\in{\cal G}$.
\end{proposition}
 \noindent{\bf Proof}. This is straightforward by the fact that 
$${\cal E}(T^n f, T^n f)\leq\frac{1}{2}\E[\sum_{i=n}^\infty\gamma_i(f,f)]$$
\fin
\begin{corollaire} Let $f\in L^2(\Omega,{\cal F}_0^\infty)$ be such that for
every $n\in N$, 
$$[0,1]^s\ni x_n\longrightarrow f(x_0,\ldots,x_n\ldots)\in\R$$
possesses a derivative in the sense of distributions in $L^2(dx_n)$ $(dx_0\cdots
dx_{n-1}dx_{n+1}\cdots)$-almost surely.

Then if 
$$\sum_{i=2}^\infty i^2(\log^\alpha i)\E[ f_i^{\prime2}]<+\infty$$
for an $\alpha>1$, then $f\in{\cal G}$.
\end{corollaire}
\noindent{\bf Proof}. The preceding proposition is here applied to the case
$(\dd,\varepsilon)=(H^1([0,1]^s,dx,\;\int \nabla^2,\;dx)$

a) Let us prove first that the hypothesis (\ref{**}) is fulfilled. For this we
use the fact that this hypothesis is satisfied on the Wiener space equipped by
the Ornstein-Uhlenbeck semi-group, as it is easily seen by the spectral
representation on the chaos. This is equivalent to saying that (\ref{**}) is
satisfied on the Gaussian structure
$$\left(\R^s,\;{\cal B}(\R^s),\;N_s(0,1),\;\int
\nabla^2,\;H^1(\R^s,N_s(0,1))\right)^{\otimes{\bf Z}}$$
with the constant $K=1$. The property is therefore true for every image
structure of this structure (cf \cite{bouleau-hirsch2}) and the result comes
then from the following easy fact :

Let be $\varphi(x)=\int_{-\infty}^x\frac{1}{\sqrt{2\pi}}e^{-\frac{y^2}{2}}dy$; 
then 
$$\int_\R(f\circ\varphi)^{\prime2}dN(0,1)\leq
\frac{1}{2\pi}\int_0^1f^{\prime2}(x)dx.$$

b) It remains only to prove that 
$$\sum_{1=2}^\infty i^2(\log^\alpha i)
a_i^2<+\infty\Rightarrow\sum_{k=0}^\infty(\sum_{i=k}^\infty a_i^2)^{1/2}<+\infty.$$
which a consequence of  the Cauchy-Schwarz inequality.\fin

\subsection{The case of  Wiener space}
\subsubsection{The Wiener space as a product space}
Let us consider the space $W=\{f\in{\cal C}([0,1],\R^d),\;f(0)=0\}$ equipped with
its Borelian $\sigma$-field ${\cal B}$ and with the Wiener measure $m$.

On the space $(\Omega,{\cal A},\PP)=\prod_{n=-\infty}^{+\infty}(W_n,{\cal
B}_n,m_n)$ where $(W_n,{\cal
B}_n,m_n)$ are copies of $(W,{\cal B},m)$ we define a Brownian motion
$(B_t)_{t\in\R_+}$ in the following manner: Letting $X_n$ be the coordinate map
from $\Omega$ into $W_n$, for $t\in]\frac{1}{2^{k+1}},\frac{1}{2^k}],\;k\in{\bf
Z}$ we put 
$$B_t=\sum_{n=k+1}^\infty\frac{X_n(1)}{2^{\frac{n+1}{2}}}+
\frac{X_k(\frac{t-1/2^{k+1}}{1/2^{k+1}})}{2^{\frac{k+1}{2}}}.$$
The process thus defined is Gaussian centred with independent increments, 
tends to zero as $t$ goes to zero and its covariance is easily computed to be
$s\wedge t$ times the identity matrix; it is therefore a standard
$\R^d$-valued Brownian motion.

The transform $\tau$ defined on $\Omega$ by 
$$X_n\circ\tau=X_{n-1}\hspace{0.5cm}n\in{\bf Z}$$
is a scaling 
$$B_t\circ\tau=\frac{1}{\sqrt{2}}B_{2t}$$
and the results of section I apply with 
$$\begin{array}{rcl}
{\cal F}_0^\infty & = & \sigma(B_s,\;s\leq 1)\\
{\cal F}_1^\infty & = & \sigma(B_s,\;s\leq \frac{1}{2})\\
{\cal F}_0^k & = & \sigma(B_s-B_{\frac{1}{2^{k+1}}},\;s\in]\frac{1}{2^{k+1}},1])
\end{array}
$$
We shall put ${\cal B}_t=\sigma(B_s,\;s<t)$.

\subsubsection{Functionals of Lipschitzian SDE's}
Let us consider maps 
$$\sigma:\R^m\times\R_+\rightarrow\R^{m\times
d},\hspace{0.5cm}b:\R^m\times\R_+\rightarrow\R^m$$
satisfying the Lipschitz hypotheses:

$\exists C>0$ such that $\forall s\in[0,1]$
$$|\sigma(x,s)-\sigma(y,s)|+|b(x,s)-b(y,s)|<C|x-y|$$
$$|\sigma(x,s)|+|b(x,s)|\leq C(1+|x|)$$
where $|\;.\;|$ is one of the equivalent norms on Euclidean spaces.

Let $X_t^x$ be the solution of the the SDE:
$$X_t^x=x+\int_0^t\sigma(X_s^x,s)dB_s+\int_0^tb(X_s^x,s)ds\hspace{0.5cm}x\in\R^m$$
\begin{proposition} Let be $f=h(X_t^x)$ for $t\leq 1$, with $h:\R^m\rightarrow
\R$ H\"{o}lderian of exponent $\lambda\in]0,1]$. Then $f\in{\cal G}.$
\end{proposition}
\noindent{\bf Proof}. Let $A$ be the H\"{o}lder constant of $h$:
$$|h(x)-h(y)|\leq A|x-y|^\lambda$$
and let $(P_t)_{t\leq 0}$ be the semi-group of the diffusion associated with
the flow $X^x_t$. By classical estimates (cf \cite{kunita} chapter 2) we have 
$$|P_u h(x)-P_u h(y)|=|\E h(X_u^x)-\E h(X_u^y)|\leq A\E|X_u^x-X_u^y|^\lambda$$
$$\leq K|x-y|^\lambda\hspace{0.5cm}\forall u\in[0,1]$$
for some constant $K$ depending on the dimensions $m,d$ and on the constants
$C$ and $A$. 

If $\varphi$ is H\"{o}lder with exponent $\lambda $, we have 
$$\mbox{var}[\varphi(X^x_s)]=\E|\varphi(X^x_s)-\E\varphi(X^x_s)|^2\leq
\E|\varphi(X^x_s)-\varphi(E(X^x_s))|^2$$
$$\leq c_1\E|X_s^x-\E X_s^x|^{2\lambda}\leq c_2(1+|x|^{2\lambda})s^\lambda$$
(cf \cite{kunita} theorem 2.1)

Now, let us remark that 
$$T^n f=\E[h(X_t^x)|{\cal
B}_\frac{1}{2^n}]\circ\tau^n=P_{t-\frac{1}{2^n}}h(X^x_{\frac{1}{2^n}})\circ\tau^n.$$
Hence by the preceding estimates we get 
$$\|T^n(f-\E
f)\|^2=\mbox{var}[T^nf]=\mbox{var}[P_{t-\frac{1}{2^n}}h(X^x_{\frac{1}{2^n}})]$$
$$\leq c(1+|x|^{2\lambda})\frac{1}{2^{n\lambda}}$$
and the series $\sum\|T^n(f-\E f)\|$ converges geometrically.\fin

\begin{proposition} Let $\mu(ds,dx)$ be a measure on $[0,1]\times \R^m$ such
that 
$$\int_{[0,1]\times \R^m}(1+|x|^\lambda)|\mu(ds,dx)|<+\infty$$
with $\lambda\in]0,1]$, and let $g$ be a H\"{o}lderian function of exponent 
$\lambda$.
Then the functional 
$$f=\int_{[0,1]\times \R^m}g(X^x_s)\mu(ds,dx)$$
belongs to ${\cal G}$.
 \end{proposition}
\noindent{\bf Proof}. We have
$$T^nf=\int_0^{\frac{1}{2^n}}\int_{\R^m}g(X^x_s)\mu(ds,dx)\circ\tau^n+\int_{\frac{1}{2^n}}^1
\int_{\R^m}P_{s-\frac{1}{2^n}}g(X_{\frac{1}{2^n}})\mu(ds,dx)\circ\tau^n$$
and hence 
$$\|T^n(f-\E f)\|\leq \int_0^\frac{1}{2^n}\int_{\R^m}\|g(X_s^x)-g(\E
X_s^x)\||d\mu|+\int_\frac{1}{2^n}^1\int_{\R^m}(var[P_{s-\frac{1}{2^n}}g(X_{\frac{1}{2^n}})])
^\frac{1}{2}|d\mu|$$
and therefore by the estimates used in the preceding proof:
$$\|T^n(f-\E f)\|\leq
\int_0^\frac{1}{2^n}\int_{\R^m}A(1+|x|^\lambda)s^\frac{\lambda}{2}|\mu(ds,dx)|
+\int_\frac{1}{2^n}^1\int_{\R^m}B(1+|x|)^\lambda)\frac{1}{2^{\frac{n\lambda}{2}}}
|\mu(ds,dx)|.$$
By hypothesis the second term is bounded by
$C\frac{1}{2^{\frac{n\lambda}{2}}}$. For the first one, let us remark that 
$$\sum_{n=0}^\infty 1_{[0,\frac{1}{2^n}]}(s)s^\frac{\lambda}{2}$$
is bounded on $s\in[0,1]$
from which it follows 
$\sum\|T^n(f-\E f)\|<+\infty$\fin

\subsubsection{Multiple Wiener integrals}
The case of multiple Wiener integrals is important on one hand because their
family is in some sense the universal diffusion process (cf \cite{azencott}
\cite{benarous}) and on the other hand because most of them are quite
irregular and such that every Borelian version is discontinuous at every point
in the Wiener space. Such functionals are not Riemann integrable and have to
be approximated by more regular functionals before simulation (cf
\cite{bouleau2}).

 Here we give some examples to illustrate which irregularity at
the origin can have functions in the Gordin class for scaling.

Let  
$$F=\int_{0<t_1<\cdots<t_m<1}h(t_1,\ldots,t_m)\;dB_{t_1}^{i_1}dB_{t_2}^{i_2}\cdots
dB_{t_m}^{i_m}$$
where $i_k\in\{1,2,\ldots,d\}$ for $k=1,\ldots,m$
with 
$$\int_{0<t_1<\cdots<t_m<1}h^2(t_1,\ldots,t_m)\;dt_1dt_2\cdots dt_m<+\infty.$$
 One has easily 
$$T^n
F=\int_{0<t_1<\cdots<t_m<1}\;\frac{1}{2^\frac{nm}{2}}h(\frac{t_1}{2^n},\ldots,
\frac{t_m}{2^n})\;dB_{t_1}^{i_1}dB_{t_2}^{i_2}\cdots
dB_{t_m}^{i_m}.$$
Therefore $F$ belongs to the Gordin class if and only if 
$$\sup_N \int_{0<t_1<\cdots<t_m<1}(\sum_{n=0}^N \frac{1}{2^\frac{nm}{2}}h(\frac{t_1}{2^n},\ldots,
\frac{t_m}{2^n}))^2 \;dt_1dt_2\cdots dt_m<+\infty$$

\noindent{\bf Example 1}. Let us take $m=1$, $h(x)=\frac{1}{x^\alpha}$,
$\alpha<\frac{1}{2}$. It is easily seen that 
$$F=\int_0^1\frac{1}{t^\alpha}\;dB_t\;\in {\cal
G}\;\;\forall\alpha<\frac{1}{2}.$$
 
\noindent{\bf Example 2}. Let us take 
$$h(x)=\frac{1}{\sqrt{x}(-\log x)^\beta}\hspace{0.5cm}\mbox{
with }\hspace{0.5cm}\beta>\frac{1}{2}.$$ Then 
$$F=\int_0^\frac{1}{2}\;\frac{1}{\sqrt{t}(-\log t)^\beta} \;dB_t$$ is in the
Gordin class if $\beta>1$, but $F\in\hspace{-0.8em}/{\cal G}$ if
$\beta\in]\frac{1}{2},1]$ although $h\in L^2[0,1]$ in that case.

\noindent{\bf Example 3}. If we take
$$h(x)=\frac{1}{\sqrt{x}}\frac{\sin(\pi\log_2 x)}{\log x}$$
the functional $F=\int_0^\frac{1}{2}h(t)\;dB_t$ gives an example of a
functional in ${\cal G}$ such that 
$$\sum_n\|T^n F\|=+\infty$$
and such that $\int_0\frac{1}{2}|h(t)|\;dB_t\in\hspace{-0.8em}/{\cal G}.$

\noindent{\bf Example 4}. Let us consider a real Brownian motion $(d=1)$, and a
function $F$ square integrable with the following Wiener chaos expansion:
$$F=F_0+\sum_m
F_m=F_0+\sum_{m=1}^\infty\int_{0<t_1<\cdots<t_m<1}h(t_1,\ldots,t_m)\;dB_{t_1}dB_{t_2}\cdots
dB_{t_m}$$
and let us suppose $|h_m(t_1,\ldots,t_m)|\leq
a_m\frac{1}{t_1^{\alpha_1^m}\cdots t_m^{\alpha_m^m}}$
with $\alpha_i^m<\frac{1}{2}$ $\forall i=1\ldots,m.$ 
We get, by the fact that the chaos are invariant by $T$, 
$$\|\sum_{n=0}^N T^n(F-F_0)\|^2=\sum_{m=1}^\infty\|\sum_{n=0}^N T^n F_m\|^2$$

$$\leq
\sum_{m=1}^\infty\frac{a_m^2}{[1-2^{(\alpha_1^m+\cdots+\alpha_m^m-\frac{m}{2})}]^2
\prod_{i=1}^m(i-2\sum_{k=1}^i\alpha_k^m)}$$
so that, if all $\alpha_i^m$'s are equal to $\alpha<\frac{1}{2}$, $F\in{\cal G}$
as soon as the series
$$\sum_{m=1}^\infty\frac{a_m^2}{m!(1-2\alpha)^m}$$
converges.
\subsection{Other factorisations of the Wiener space.}
\subsubsection{}
Let $(\chi_n(t))_{n\geq 0}$ be an orthonormal basis of $L^2[0,1]$ and let
$\varphi_n(t)=\int_0^t \chi_n(s)\;ds$. Let $(g_n)_{n\geq 0}$ be a sequence of
independent standard Gaussian variables built as the coordinates of 
$(\Omega,{\cal A}, \PP)=(\R^\N,{\cal B}(\R^\N),N(0,1)^{\otimes\N})$.

The series 
\begin{equation}
\label{31}
\sum_n g_n\varphi_n
\end{equation}
converges in ${\cal C}([0,1])$ a.s. and in $L^p((\Omega,{\cal A}, \PP),{\cal
C}([0,1]))$ $p\in[1,\infty[$ and its sum is a Brownian motion under $\PP$.

Indeed, if on the Wiener space we put 
$\tilde{\chi}_n(\omega)=\int_0^1 \chi_n(s)\;dB_s$
and ${\cal F}_n=\sigma(\tilde{\chi}_k,\;k\leq n)$, we obtain, denoting by $B$
the identity map from ${\cal C}[0,1]$ into itself,
\begin{equation}\label{truc}
\E[B|{\cal F}_n]=\sum_{k=0}^n\tilde{\chi}_k \varphi_k 
\end{equation}
as can be seen by applying a continuous linear functional $\mu$ on ${\cal
C}[0,1]$ to both sides of (\ref{truc}) and by remarking that $(\mu,
\tilde{\chi}_0,\ldots,\tilde{\chi}_n)$ is a Gaussian array. By the convergence
properties of vector martingales, we have therefore
\begin{equation}\label{32}
B=\sum_{k=0}^\infty\tilde{\chi}_k\varphi_k
\end{equation}
 a.s. and in $L^p$. Since the family
of partial sums of the series (\ref{31}) has the same law as the sums of
(\ref{32}) the assertion is proved.

Such a representation of the Brownian motion 
$$B=\sum_{k=0}^\infty\tilde{\chi}_k\varphi_k$$
allows us to define the shift, and the associated Gordin class clearly depends 
on the basis $(\chi_n)$ which is chosen.

\subsubsection{}
 The case of Haar functions is particularly interesting. Let us put
$\chi=1_{[0,\frac{1}{2}[}-1_{[\frac{1}{2},1[}$
and 
\begin{eqnarray}
\chi_{m,k}(t) & = & 2^{\frac{m}{2}}\chi(2^mt-k)\\
\varphi_{m,k}(t) & = & \int_0^t\chi_{m,k}(s)\;ds
\end{eqnarray}
for $t\in\R_+$, $m\in{\bf Z}$, $k\in\N$.

The functions $(\chi_{m,k})_{m\in{\bf Z},\;k\in\N}$ form an orthonormal basis
of $L^2(\R_+)$ and if $g_{m,k}$ are standard independent Gaussian variables,
the Brownian motion can be represented by 
$$B(\omega,t)=\sum_{m=-\infty}^{+\infty}(\sum_{k=0}^\infty
\varphi_{m,k}(t)g_{m,k}(\omega))$$
and the scaling  studied in  paragraph II.3 is the mapping which transforms
the sequence
$$(g_{m,k}(\omega))_{m,k}$$
into the sequence
$$(g_{m-1,k}(\omega))_{m,k}.$$
The space generated by the functions $(\chi_{m,k})_{m\geq 0,\;0\leq k<2^m}$ is
the subspace of $L^2[0,1]$ orthogonal to the constants, and the process 
\begin{equation}\label{33}
Z_t=\sum_{m=0}^\infty\sum_{k=0}^{2^m-1}\varphi_{m,k}(t) g_{m,k}
\end{equation}
is a standard Brownian bridge vanishing at zero and one. The representation
(\ref{33}) is unique and converges in ${\cal C}_{00}[0,1]=\{f\in{\cal
C}[0,1]\;f(0)=f(1)=0\}$. The functions $\varphi_{m,k}$ form a Schauder basis
of this space. If $f\in {\cal C}_{00}[0,1]$ with 
\begin{equation}\label{57}
f(t)=\sum_{m=0}^\infty\sum_{h=0}^{2^m-1}\varphi_{m,k}(t)a_{m,k}(f)
\end{equation}
there holds
$$a_{m,k}(f)=[2f(\frac{k}{2^m}+\frac{1}{2^{m+1}})-f(\frac{k}{2^m})-f(\frac{k+1}{2^m})]
2^{\frac{m}{2}}.$$
The Banach spaces of H\"{o}lderian functions of exponent $\alpha\in]0,1[$ of 
${\cal C}_{00}[0,1]$ can be interpreted in terms of spaces $\ell^\infty$ and
$c_0$ on the sequences $(2^{m\alpha}a_{m,k})_{m,k}$ (cf \cite{ciesielski}).

To approach a continuous function by a partial sum of the series (\ref{57}) is
convenient practically, and if we change the notations by putting 
$a_{2^m+k}=a_{m,k}$ $m\geq 0$, $k=0,\ldots,2^m-1$ the simple shift on the $a_n$
i.e., the transform
$$F(a_1,\ldots,a_n,\ldots)\longrightarrow
F\circ\tau=F(a_0,a_1,\ldots,a_{n+1},\ldots)$$
(which does not correspond to a scaling) is quite thrifty in random drawings.
By proposition 11, a sufficient condition for a function $F$ to be in the
Gordin class for this transform is that it be in $L^2$ and possess partial
derivatives such that 
$$\sum_{k=0}^\infty(\sum_{i=k}^\infty\E(F_i^{\prime 2}))^\frac{1}{2}<+\infty$$
where the expectation is taken on $(\R^\N,{\cal B}(\R^\N),N(0,1)^{\otimes\N})$.

\noindent{\bf Example}. For fixed $t\in[0,1]$, let us consider the functional
$$F(\omega)=\sum_{n\geq
0}\frac{1}{n+1}\sqrt{\varphi_n(t)}\left(\int_0^1\chi_n(s)\;dZ_s\right)^2$$
which, with the preceding notations $F$ can be written
$$F=\sum_{n\geq 0}\frac{1}{n+1}\sqrt{\varphi_n(t)}\;a_n^2.$$
Now $F$ belongs to $L^2$ by the fact that the series 
$$\sum_{n\geq 0}\frac{1}{(n+1)^2}\varphi_n(t)$$
converges and we have $F_i^\prime =\frac{2}{i+1}\sqrt{\varphi_i(t)}a_i$ so that 
$$\sum_{k=0}^\infty\sum_{i=k}^\infty\E F_i^{\prime
2}=\sum_{k=0}^\infty\frac{4}{k+1}\varphi_k(t)<+\infty$$
because $(\frac{1}{k+1})\in \ell^2$. And thus $F\in{\cal G}$.


\begin{thebibliography}{40}
{\small

\bibitem[1] {benalaya}M.B. Alaya. On the simulation of expectations of random variables depending on a
stopping time, {\it Stoc. Anal. and Appl.} 11 (1993) 133-153

 \bibitem[2] {azencott}R. Azencott. Formule de Taylor
stochastique et d\'{e}veloppement asymptotique d'int\'{e}grales de Feynmann.\\
 237-285, in {\it Sem.
Prob. XVI supp. G\'{e}om\'{e}trie diff. stoch.} Lect. Notes in M. 921 Springer
(1982)

\bibitem[3] {benarous}G. Ben Arous. Flots et s\'{e}ries de Taylor
stochastiques.\\
{\it Prob. Th. Rel. Fields,} 81, 29-77, (1989)

\bibitem[4] {berger}E. Berger. An almost sure invariance principle for
stationary ergodic sequences of Banach space valued random variables.\\
{\it Prob. Th. Rel. Fields.} 84, 161-201, (1990)
 
\bibitem[5] {bouleau}N. Bouleau. On effective computation of expectations
in large or infinite dimension.\\
{\it J. of Computational and App. Math.} 31,23-34, (1990)

\bibitem[6]{bouleau2}N. Bouleau. Irregular and simulatable functionals on Wiener space, in {\it
Probabilit\'{e}s Num\'{e}riques}, p39-53, INRIA, (1991)

 \bibitem[7] {bouleau-hirsch1}N. Bouleau, F.
Hirsch. Formes de Dirichlet g\'{e}n\'{e}rales et densit\'{e} de variables al\'{e}atoires sur l'espace de
Wiener.\\
 {\it J. Funct. Analysis,} vol 69, 227-259, (1986)

\bibitem[8] {bouleau-hirsch2}N. Bouleau, F. Hirsch. Alg\`{e}bre des
structures de Dirichlet.\\
 {\it C. R. Acad. Sc. Paris} t310, sI, 15-18, (1990)

\bibitem[9] {bouleau-pages}N. Bouleau, G. Pag\`{e}s, J. Xiao. Extension des
M\'{e}thodes de Monte Carlo acc\'{e}l\'{e}r\'{e}es.\\
Contrat DRET, CERMA Ecole Nationale des Ponts et Chauss\'{e}es, Paris (1990)
 
\bibitem[10]{ciesielski}Z. Ciesielski. On the isomorphism of spaces $H_\alpha$ and $m$.\\
{\it Bull. Acad. Polonaise des Sc., s. des Sc. math. astr. phys.}, VIII, n¡4,  217-222, (1960)



\bibitem[11] {cornfeld}I.P. Cornfeld, S.V. Fomin, Ya.G. Sinai. {\it Ergodic
Theory}.\\
Springer (1982)


\bibitem[12] {h.faure}H. Faure.Discr\'{e}pance de suites associ\'{e}es \`{a}
un syst\`{e}me de num\'{e}ration (en dimension s).\\ {\it Acta Arithm.} 41,
337-351, (1982)
 


\bibitem[13] {gordin}M.I. Gordin. The central limit theorem for stationary
sequences.\\
 {\it Soviet Math. Dokl.} vol. 10, n¡5, 1174-1175,  (1969)

\bibitem[14] {halasz}G. Hal\'{a}sz. Remarks on the remainder in Birkhoff's
ergodic theorem.\\
{\it Acta Math. Acad. Hungar.} 28,389-395, (1976)

\bibitem[15] {hall-heyde}P. Hall, C.C. Heyde. {\it Martingale limit theory and
applications}.\\ Acad. Press (1980)
 
\bibitem[16] {heyde-scott}C.C. Heyde, D.J. Scott. Invariance principles for
the law of iterated logarithm for martingales and processes with
stationary increments.\\
{\it Ann. Prob.} vol 1, n¡3, 428-436, (1973)

\bibitem[17] {krengel1}U. Krengel. On the speed of convergence in the
ergodic theorem.\\
{\it Monatsh. Math.} 86, 3-6, (1978)

\bibitem[18] {krengel2}U. Krengel. {\it Ergodic theorems}.\\
 de Gruyter (1985)

\bibitem[19] {kunita}H. Kunita. Stochastic differential equations and
stochastic flows of diffeomorphisms.\\
in {\it Ecole d'\'{e}t\'{e} de St Flour XII,} 143-303, Springer lect. notes in M.
1087, (1982)
 

\bibitem[20] {niederreiter1}H. Niederreiter.Quasi-Monte Carlo methods and
pseudo-random numbers.\\
{\it Bull. Amer. Math. Soc.} 84, 957-1041, (1978)
 
\bibitem[21] {niederreiter2}H. Niederreiter.Low-discrepancy and
low-dispersion sequences.\\{\it J. of Number Th.} 30, n¡1, 51-70, (1988)

 
\bibitem[22] {sarkar}P.K. Sarkar, M.A. Prasad. A comparative study of
pseudo and quasi random sequences for solution of integral equations .
{\it J. of computational Physics} 68, 66-88, (1987)
 
\bibitem[23] {scott}D.J. Scott. Central limit theorems for martingales and
processes with stationary increments using a Skorokod representation
approach.
{\it Adv. Appl. Prob.} 5,119-137, (1973)


  }
\end{thebibliography}
 \end{document}